\documentclass[journal,twoside]{IEEEtran}
\usepackage{color}
\usepackage{cite}
\ifCLASSINFOpdf
   \usepackage[pdftex]{graphicx}
\else
\fi
\usepackage[cmex10]{amsmath}
\usepackage{amssymb}
\usepackage{amsfonts}
\usepackage{float}
\usepackage{subfigure}
\usepackage{multirow}
\usepackage[capitalize]{cleveref}
\usepackage{algorithm}
\usepackage{algpseudocode}
\usepackage{multirow}
\usepackage{multicol}
\usepackage{todonotes}

\usepackage{tikz}
\tikzstyle{cnode} = [draw, circle,scale=.6]
\tikzstyle{level 1} = [level distance=.35\textwidth, sibling distance=1\textwidth]
\tikzstyle{level 2} = [level distance=.35\textwidth, sibling distance=.45\textwidth]
\tikzstyle{level 3} = [level distance=.3\textwidth, sibling distance=.25\textwidth]
\usetikzlibrary{decorations.pathreplacing,patterns,external}

\usetikzlibrary{external}
\newcommand{\rev}[1]{{\color{black}{#1}}}
\newcommand{\newrev}[1]{{\color{black}{#1}}}

\ifCLASSOPTIONcompsoc
 \usepackage[caption=false,font=normalsize,labelfont=sf,textfont=sf]{subfig}
 \else
 \usepackage[caption=false,font=footnotesize]{subfig}
 \fi

\begin{document}
\title{A Butterfly-Accelerated Volume Integral Equation Solver for Broad Permittivity and Large-Scale Electromagnetic Analysis}
\author{Sadeed~B.~Sayed, \textit{Member}, \textit{IEEE}, Yang Liu, \textit{Member}, \textit{IEEE}, Luis J. Gomez, \textit{Member}, \textit{IEEE}, and Abdulkadir C. Yucel, \textit{Senior Member}, \textit{IEEE}
\thanks{Manuscript received March 3, 2021. This work was supported by Nanyang Technological University under a Start-Up Grant. (\textit{Corresponding author: Abdulkadir C. Yucel.})}
\thanks{S.~B.~Sayed and A. C. Yucel are with the School of Electrical and Electronic Engineering, Nanyang Technological University, Singapore 639798. (e-mails: sadeed.sayed@ntu.edu.sg, acyucel@ntu.edu.sg).}

\thanks{Y. Liu is with Computational Research Division, Lawrence Berkeley National Laboratory, Berkeley, CA 94720, USA. (e-mail: liuyangzhuan@lbl.gov). }
\thanks{L. J. Gomez is with the division of Electrical and Computer Engineering, Purdue University, USA. (e-mail: ljgomez@purdue.edu).}
}
\markboth{IEEE TRANSACTIONS ON ANTENNAS AND PROPAGATION,~Vol.~P, No.~PP, Month~Year}
{SAYED \MakeLowercase{\textit{et al.}}: A Butterfly-Accelerated VIE Solver }

\maketitle

\begin{abstract}
A butterfly-accelerated volume integral equation (VIE) solver  is proposed for fast and accurate electromagnetic (EM) analysis of scattering from heterogeneous objects. The proposed solver leverages the hierarchical off-diagonal butterfly (HOD-BF) scheme to construct the system matrix and obtain its approximate inverse, used as a preconditioner. Complexity analysis and numerical experiments validate the $O(N\log^2N)$ construction cost of the HOD-BF-compressed system matrix and $O(N^{1.5}\log N)$  inversion cost for the preconditioner, where $N$ is the number of unknowns in the high-frequency EM scattering problem. For many practical scenarios, the proposed VIE solver requires less memory and computational time to construct the system matrix and obtain its approximate inverse compared to a $\mathcal{H}$ matrix-accelerated VIE solver. The accuracy and efficiency of the proposed solver have been demonstrated via its application to the EM analysis of large-scale canonical and real-world structures comprising of broad permittivity values and involving millions of unknowns.  
\end{abstract}


%
\IEEEpeerreviewmaketitle

\section{Introduction}\label{sec:intro}
%
%
%
%

\IEEEPARstart{T}{he} analysis of electromagnetic (EM) scattering from  heterogeneous scatterers is traditionally performed via volume integral equation (VIE) solvers \cite{schaubert1984tetrahedral,Sancer_2006,Botha_2006}. These solvers discretize the VIE and obtain a dense system of linear equations with $N$ unknowns. Oftentimes, these solvers' applicability to the large-scale problems is limited due to  $O(N^2)$ or  $O(N^3)$ memory and CPU time requirements for the iterative or direct solution of the dense VIE system, respectively. On the other hand, differential equation solvers require the solution of a sparse system of linear equations and consequently demand less computational resources compared to VIE solvers when applied to the large-scale problems. That said, unlike VIE solvers, differential equation solvers do not implicitly satisfy the radiation conditions and can suffer from the numerical dispersion. 

So far, a plethora of methods has been developed to lower the memory and CPU time requirements of VIE solvers \cite{Sertel2004FMMVIE,Lu2003FMMVIE,jacob_precorrectedFFT,Nie2005P-FFT,yucel2018voxhenry,Ozdemir2004LRVIE,Wenwen_Chai_2009,chai2010,chen2019tensor}. These methods were primarily developed to expedite the matrix-vector multiplications during the iterative solution of the VIE system. They leverage  analytic or algebraic compression techniques, including the fast multipole method (FMM) \cite{Sertel2004FMMVIE,Lu2003FMMVIE}, fast Fourier transform \cite{jacob_precorrectedFFT,Nie2005P-FFT,yucel2018voxhenry,wang2020voxcap}, low-rank compression \cite{Ozdemir2004LRVIE}, hierarchical matrices as kernel-free FMM \cite{Wenwen_Chai_2009,chai2010}, and, more recently, tensor decompositions \cite{chen2019tensor}. Despite the success of these methods in expediting matrix-vector multiplications during iterative solution, the discretized VIE systems are often ill-conditioned. This results in a prohibitively large iteration count during the iterative solutions of the VIE systems. In particular, the systems are ill-conditioned in scenarios where the scatterers have large permittivity  \cite{markkanen2012analysis,Costabel_2012,Kottmann_2000,gomez2017internally,gomez2015low}, negative permittivity  \cite{budko2006spectrum,liu2020fastskeletonization,gomez2015volume,sayed2020multi}, or large electrical size. This results in several challenges for using VIE in the EM analysis of biomedical procedures applied to high-permittivity tissues \cite{gomez2017icvsie}, EM analysis of the structures covered by negative permittivity plasma \cite{yucel2018internally}, and EM analysis of electrically-large photonic devices \cite{groth2019circulant}. For these scenarios, fast matrix-vector multiplication schemes have to be used in conjunction with effective preconditioners. Alternatively, fast direct methods can be utilized.

Existing fast direct solvers for the VIE include Huygens' equivalence principle-based algorithms \cite{Chew1993Nepal,Chew1993Huygens} yielding an $O(N^2)$ complexity, and low-rank-based algorithms such as $\mathcal{H}$ matrices \cite{gholami2018h, gholami2020surface}, $\mathcal{H}^2$ matrices \cite{omar2015linear,ma2017accuracy,ma2019direct}, and skeletonization \cite{greengard_gueyffier_martinsson_rokhlin_2009,Corona2015ON,liu2020fastskeletonization}. These low-rank-based algorithms can attain $O(N)$ or $O(N\log N)$ complexities for static or low-frequency scattering problems, but tend to become less efficient as the electrical size of the scatterer increases. Low-rank-based direct solvers have also been developed for surface integral equations (SIE) \cite{hackbusch2003introduction,Wenwen_Chai_2009,chai2010}, but low-rank algorithms tend to perform better for VIE solvers. This is due to the fact that the density of discretization elements per wavelength is typically much larger for the VIE compared to that for the SIE.  

In contrast to the low-rank algorithms, here we consider another class of algebraic compression techniques called butterfly algorithms \cite{Eric_1994_butterfly,michielssen_multilevel_1996,li2015butterfly,Yingzhou_2017_IBF,Pang2020IDBF} for fast iterative and/or direct solution of the VIE. For a $m\times n$ matrix, its butterfly representation can be treated as an extension of the fast Fourier transform into an $O(\log n)$-factor multiplicative matrix decomposition that is well-suited to compress highly-oscillatory operators \cite{Candes_butterfly_2009,Oneil_2010_specialfunction,bremer2020SHT}. Each butterfly factor contains $O(n)$ number of small dense blocks whose sizes are oftentimes referred to as butterfly ranks. When combined with hierarchical matrix techniques, they have been leveraged for the direct or iterative solutions of  2-dimensional \cite{michielssen_multilevel_1996,Han_2013_butterflyLU,Liu_2017_HODBF,Yang_2020_BFprecondition}, 3-dimensional \cite{Han_2017_butterflyLUPEC,Han_2018_butterflyLUDielectric,Yaniv2021BF} SIE, and differential equations \cite{liu2020MF} for electrically-large scattering problems. In general, the butterfly-accelerated SIE solvers can be classified into the strong-admissibility condition-based ones (such as hierarchical LU \cite{Han_2017_butterflyLUPEC}) and weak-admissibility conditioned-based ones (such as hierarchically off-diagonal butterfly (HOD-BF) algorithm \cite{Liu_2017_HODBF}). The former compresses only matrix blocks for the well-separated interactions, while the latter compresses the matrix blocks for both nearby and well-separated interactions  \cite{hackbusch2004hierarchical}. \rev{They can be treated as the butterfly extensions of the strong-admissibility $\mathcal{H}$  \cite{gholami2018h} and weak-admissibility hierarchically off-diagonal low-rank (HOD-LR) \cite{ambikasaran2013mathcal,AMINFAR2016170,Shaeffer2018hodlr} matrices, respectively.}   

This paper proposes a butterfly-accelerated VIE solver for the EM scattering analyses of large-scale scatterers comprising broad permittivity values. The proposed solver leverages the HOD-BF scheme to compress the blocks in the VIE system and obtain the approximate inverse of the system matrix for forming a preconditioner. During the iterative solution of the VIE system, the matrix-vector multiplications are accelerated using the HOD-BF-compressed blocks while the iteration count is dramatically dropped to single digits with the preconditioner, even for highly ill-conditioned problems.  Compared to the strong-admissibility solvers \cite{Han_2017_butterflyLUPEC,Han_2018_butterflyLUDielectric}, the proposed solver has simpler butterfly arithmetic, smaller leading constants in complexity, and significantly better parallelization performance. The HOD-BF scheme leveraged in the proposed solver was first introduced \rev{as a sequential algorithm} to solve 2-dimensional SIE with $O(N\log^2N)$ memory and $O(N^{1.5}\log N)$ CPU time requirements \cite{Liu_2017_HODBF}; the same complexity can also be attained for 3-dimensional SIE despite of the non-uniform rank patterns due to weak-admissibility \cite{liu2020butterfly}. To the best of our knowledge, neither the HOD-BF scheme nor other butterfly schemes has been applied to the 3-dimensional VIE before, and consequently the complexity analysis of HOD-BF scheme has never been performed for the VIE. \rev{Unlike the HOD-LR scheme that scales poorly for high-dimensional problems, e.g., $O(N^{5/3})$ memory and $O(N^{7/3})$ CPU for 3D VIEs, we show that HOD-BF scheme's complexity is dimension independent.} Furthermore, the advantages of the HOD-BF scheme compared to the  low-rank-based $\mathcal{H}$ matrix scheme have never been shown. \rev{All the numerical experiments are carried out  using a distributed-memory parallel implementation of HOD-BF, briefly explained in this study. Without such implementation, the demonstration of the advantages of HOD-BF compared to the HOD-LR and $\mathcal{H}$ matrix schemes would not be possible.} These are the novelties introduced in this paper. 

Here, we theoretically and numerically analyze the rank scaling and complexity of the proposed HOD-BF-accelerated VIE solver. Furthermore, we extensively compare the performance of the proposed VIE solver with that of a low-rank-based $\mathcal{H}$ matrix-accelerated VIE solver.  The numerical analysis shows that the proposed solver achieves $O(N\log^2N)$ memory and CPU time scaling for constructing the HOD-BF-compressed system matrix as well as $O(N\log^2N)$  memory and  $O(N^{1.5}\log N)$ CPU time scaling for inverting the system matrix to obtain the preconditioner. Compared to the  $\mathcal{H}$ matrix-accelerated VIE solver, the proposed solver requires 1.4/3.5 times less CPU time and 3.8/3.9 times less memory for the construction/inversion when used for the EM analysis of a dielectric sphere involving $0.9$ M unknowns. Similar comparison performed on a NASA almond involving $1.6$ M unknowns shows that the proposed simulator requires 5.0/3.2 times less CPU time and 3.8/2.7 times less memory for the construction/inversion compared to the  $\mathcal{H}$ matrix-accelerated VIE solver. Moreover, the proposed solver with the preconditioner requires maximum five iterations during the iterative solution of an highly ill-conditioned VIE system generated for the EM scattering analysis of negative permittivity sphere involving $1.7$ M unknowns. 

The rest of the paper is organized as follows: \cref{sec:form} describes the VIE formulation. \cref{sec:HODBF} presents the construction of the HOD-BF-compressed system matrix via butterfly algorithm, inversion of the compressed system matrix to obtain preconditioner, complexity estimates, and parallelization strategy. \cref{sec:num_res} provides numerical results to demonstrate the scalability, accuracy, and efficiency of the proposed butterfly-accelerated VIE solver. The conclusions are drawn in Section \ref{sec:con}. Throughout this paper, a time-harmonic variation of field quantities is assumed and suppressed.

\section{VIE Formulation}\label{sec:form}

Consider a time-harmonic EM field with incident electric and magnetic fields,  ${\mathbf{E}}^{\mathrm{inc}}(\mathbf{r})$ and  $\mathbf{H}^{\mathrm{inc}}({\mathbf{r}})$, at frequency $f$. The field impinges on a heterogeneous scatterer ${V}$ residing in a background medium with permittivity $\varepsilon_0$,  permeability $\mu_0$, and intrinsic impedance ${\eta_0}=\sqrt{\mu_0/\varepsilon_0}$. Here the background medium is assumed to be free-space. Using the volume equivalence principle, the  VIE for the electric flux density ${\mathbf{D}}({\mathbf{r}})$ is obtained as \cite{schaubert1984tetrahedral}
\begin{align}
 {\mathbf{D}}({\mathbf{r}})/\varepsilon ({\mathbf{r}}) -  j \omega\eta_0\mathcal{L}\left[ \kappa{{\mathbf{D}}} \right](\mathbf{r}) = {{\mathbf{E}}^{{\mathrm{inc}}}}({\mathbf{r}});{\mathbf{r}} \in V,
\label{eq:1}
\end{align}
where 
\begin{align}
 \mathcal{L}[{\mathbf{X}}]({\mathbf{r}})=-j{k_0}\int_V dv' \left[ \mathbf{I}+\frac{1}{k_0^2}\nabla\nabla\cdot
 \right] {\mathbf{X}}({\mathbf{r}}) {\mathbf{g}}(\mathbf{r},\mathbf{r}').
\label{eq:2}
\end{align}
Here $\omega=2\pi f$ is the angular frequency, $\mathbf{r}$ and $\mathbf{r'}$ are  the observation and source points, $k_0=\omega\sqrt{\mu_0\varepsilon_0}=2\pi/\lambda_0$ is the wave number, $\lambda_0$ is the wavelength, $\varepsilon(\mathbf{r})$ is the dielectric permittivity, $\kappa(\mathbf{r})=\left(\varepsilon(\mathbf{r})-\varepsilon_0\right)/\varepsilon(\mathbf{r})$ is the dielectric contrast, and \rev{ ${\mathbf{g}}(\mathbf{r},\mathbf{r}')={e^{-jk_0|\mathbf{r}-\mathbf{r}'|}}/{4\pi|\mathbf{r}-\mathbf{r}'|}$} is the Green's function. 

To solve the VIE, the scatterer is discretized by a volumetric tetrahedral mesh having ${N}$ faces. Then, the electric flux in ${V}$ is approximated via Schaubert Wilton Glisson (SWG) basis functions ${{\mathbf{f}_n(\mathbf{r})}}$ as \cite{schaubert1984tetrahedral}
\begin{align}
 {\mathbf{D}}({\mathbf{r}}) = \sum\limits_{n{\mathrm{ = 1}}}^{{N}} {D_n{\mathbf{f}}_n}({\mathbf{r}});{\mathbf{r}} \in V.
\label{eq:3}
\end{align}
Note that one SWG is defined for each of the ${N}$ faces of the tetrahedral mesh, and each basis function ${{\mathbf{f}_n(\mathbf{r})}}$ has support on $V_n$, the volume of the tetrahedron which the face $N$ belongs to. 
Substituting  (\ref{eq:3}) into (\ref{eq:1}) and applying Galerkin testing yields the discretized VIE linear system of equations as
\begin{align}
 {\mathcal{Z}\mathcal{I}} = {\mathcal{V}}, 
\label{eq:4}
\end{align}
where
\begin{align}
\nonumber {{\mathcal{Z}}(m,n)} &= \left\langle {{{\mathbf{f}}_m({\mathbf{r}})},\frac{{{\mathbf{f}_n({\mathbf{r}})}}}{{{\varepsilon(\mathbf{r})}}} -  j \omega\eta\mathcal{L}\left(\kappa(\mathbf{r}) {{{\mathbf{f}}_n({\mathbf{r}})}} \right)} \right\rangle_m  \\  
\nonumber {{\mathcal{I}}(n)} &= D_n \\
 {{\mathcal{V}}(m)} &= \left\langle {{{\mathbf{f}}_m({\mathbf{r}})},{{\mathbf{E}}^{{\mathbf{inc}({\mathbf{r}})}}}} \right\rangle_m;\hspace{6mm}n,m=1,...,N.
\label{eq:5}
\end{align}
 Here, $\mathcal{Z}$ is the system matrix, $\mathcal{I}$  is the unknown coefficient vector, and $\mathcal{V}$  is the excitation vector. Furthermore,  $\kappa(\mathbf{r}) = (\varepsilon(\mathbf{r})-\varepsilon_0)/\varepsilon(\mathbf{r})$ and   $<\mathbf{a},\mathbf{b}>_m=\int_{V_m} \mathbf{a}\cdot\mathbf{b} d\mathbf{r}$ denotes the standard inner product of $\mathbf{a}$ with $\mathbf{b}$ in $V_m$.  
 
\newrev{ 
The linear system \eqref{eq:4} is solved with preconditioned iterative solvers as  
 \begin{align}
\mathcal{Z}^{-1}_{\chi_{fact}}{\mathcal{Z}_{\chi_{con}}\mathcal{I}} = \mathcal{Z}^{-1}_{\chi_{fact}}{\mathcal{V}}. 
\label{eq:precon}
\end{align} 
Here $\mathcal{Z}_{\chi_{con}}$ is a compressed representation of $\mathcal{Z}$ with tolerance $\chi_{con}$ using e.g., $\mathcal{H}$ \cite{gholami2018h}, HOD-LR \cite{ambikasaran2013mathcal}, or the proposed HOD-BF \cite{Liu_2017_HODBF} techniques. $\mathcal{Z}^{-1}_{\chi_{fact}}$ is the compressed representation of $\mathcal{Z}^{-1}$ with tolerance ${\chi_{fact}}\geq{\chi_{con}}$. We will drop the subscripts of $\mathcal{Z}_{\chi_{con}}$ and $\mathcal{Z}^{-1}_{\chi_{fact}}$ in the rest of this paper. The preconditioned linear system \eqref{eq:precon} can be solved via an iterative solver (such as the transpose-free quasi-minimal residual method (TFQMR) \cite{Roland1993TMQMR}) by setting a stopping criteria ${\chi_{sol}}$. 
}

\section{HOD-BF Acceleration and Preconditioner\label{sec:HODBF}}
\newrev{The HOD-BF scheme first generates a hierarchical partitioning of the system matrix and then obtains the compressed off-diagonal blocks representing non-overlapping interactions as the butterfly format. The compressed blocks permit rapid matrix-vector multiplications and inversion, particularly for high-frequency regimes and 3-dimensional scatterers.} In what follows, we briefly summarize the construction of the system matrix via HOD-BF-compressed blocks and inversion of the compressed system matrix to obtain preconditioner. 
\subsection{Construction of HOD-BF-compressed system matrix}
\begin{figure}[!htpb]
	\centering
	\includegraphics[width=0.95\columnwidth]{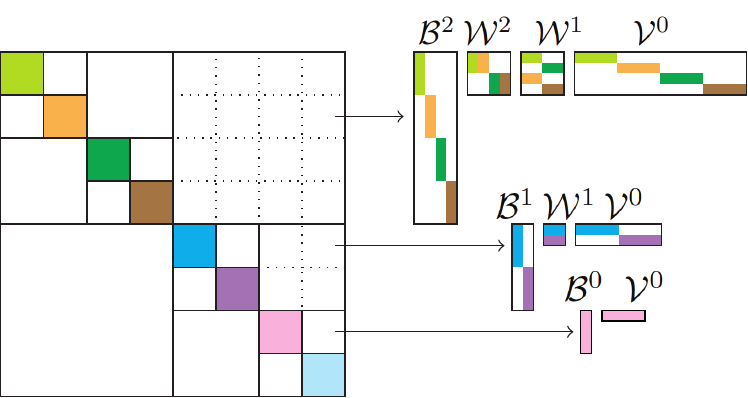}
 \caption{	
	Example of a $L_H=3$-level HOD-BF-compressed matrix. At the leaf level, the diagonal blocks $D_\tau$ are stored as dense blocks; at levels $l>L_H$, the off-diagonal blocks are compressed as $L=L_H-l$-level butterfly representations. \newrev{The HOD-BF matrix is distributed over 8 MPI processes each denoted by one color. For each $\tau$, $D_\tau$ and $\mathcal{B}_{\tau}$ are assigned to the same set of processes. $\mathcal{B}_{\tau}$ is distributed over multiple processes using parallel layouts of \cite{liu2020butterfly}.} \label{fig:HOD-BF}}
\end{figure}

The construction begins by recursively subdividing the scatterer into two sub-scatterers with approximately equal numbers of unknowns using a tree clustering algorithm (such as KD) \cite{Han_2013_butterflyLU}, until the subscatterer contains a prescribed number of unknowns. This procedure results in a complete binary tree $\mathfrak{T}_H$ of $L_H$ levels with root level $0$ and leaf level $L_H$. Each node $\tau$ at level $l$ is an index set $\tau \subset \{1,\dots,N\}$ associated with the corresponding subscatterer. For a non-leaf node $\tau$ at level $l$ with children $\tau_1$ and $\tau_2$, $\tau=\tau_1\cup\tau_2$ and $\tau_1 \cap \tau_2 = \emptyset$. For a non-root node $\tau$, its parent is denoted $p_\tau$. \newrev{Let $\mathcal{D}_\tau = \mathcal{Z}(\tau, \tau)$. For the root $\tau$, $\mathcal{D}_\tau=\mathcal{Z}$;} for each leaf node $\tau$, $\mathcal{D}_\tau$ is directly computed as dense blocks; off-diagonal blocks are compressed using the butterfly representation described below. As an example, \cref{fig:HOD-BF} shows a HOD-BF-compressed matrix with $L_H=3$ levels.  

Specifically, let $\tau_1$ and $\tau_2$ be two siblings in $\mathfrak{T}_H$ on level $l$. These two sibling nodes correspond to two off-diagonal blocks $\mathcal{B}_{\tau_1} = \mathcal{Z}(\tau_1, \tau_2)$ and $\mathcal{B}_{\tau_2} = \mathcal{Z}(\tau_2, \tau_1)$. As an example consider the butterfly compression of the $m\times n$ block $\mathcal{B}=\mathcal{B}_{\tau_1}$ with $o=\tau_1$ and $s=\tau_2$, then $\mathcal{B}=\mathcal{Z}(o,s)$ is compressed as a butterfly with $L=L_H-l$ levels. Let $\mathfrak{T}_{o}$ and $\mathfrak{T}_{s}$ denote subtrees of $\mathfrak{T}_H$ with $L$ levels, rooted at nodes $o$ and $s$ respectively. The butterfly compression requires the {\it complementary low-rank property}: for any level $0\leq l \leq L$, any node $\tau$ at level $l$ of $\mathfrak{T}_{o}$ and any node $\nu$ at level $L-l$ of $\mathfrak{T}_{s}$, the subblock $\mathcal{Z}(\tau,\nu)$ is numerically low-rank with rank $r_{\tau,\nu}$ bounded by a small number $r$ called the butterfly rank. As will be discussed in \cref{sec:CC}, we do not necessarily require constant butterfly ranks to attain a quasi-linear representation. Given the complementary low-rank property, we can compress any subblock $\mathcal{Z}(\tau,\nu)$ above as a low-rank product in an interpolative form: 
\begin{equation}\label{eqn:ID}
\mathcal{Z}(\tau,\nu) \approx \mathcal{Z}(\tau, \bar{\nu}) \mathcal{V}_{\tau,\nu}.
\end{equation}
where $\bar{\nu}$ represents the skeleton columns, $\mathcal{Z}(\tau, \bar{\nu})$ is the skeleton matrix, and $\mathcal{V}_{\tau,\nu}$ is the interpolation matrix. The approximation in \cref{eqn:ID} is typically computed with a compression tolerance ${\chi_{con}}$ such that $|\mathcal{Z}(\tau,\nu) - \mathcal{Z}(\tau, \bar{\nu}) \mathcal{V}_{\tau,\nu}|_F=O(\chi_{con})|\mathcal{Z}(\tau,\nu)|_F$. This low-rank product can be computed from the interpolative decomposition (ID) \cite{Cheng2005OnTC}. Moreover, the interpolation matrices $\mathcal{V}_{\tau,\nu}$ for non-leaf nodes $\nu$ can be defined in a nested fashion as 
 
\begin{equation}\label{eqn:nested_basis}
\mathcal{V}_{\tau,\nu} =
\mathcal{W}_{\tau,\nu}\begin{bmatrix}
\mathcal{V}_{p_\tau, \nu_1} & \\
& \mathcal{V}_{p_\tau, \nu_2}
\end{bmatrix}.
\end{equation}
where $\mathcal{W}_{\tau,\nu}$ are referred to as the transfer matrices. \newrev{$\mathcal{W}_{\tau,\nu}$ is computed from ID: $\mathcal{Z}(\tau,\bar{\nu})\mathcal{W}_{\tau,\nu}\approx[\mathcal{Z}(\tau,\bar{\nu_1}), \mathcal{Z}(\tau,\bar{\nu_2})]$}

Once all the interpolation and transfer matrices are computed, the butterfly factorization of $\mathcal{B}$ can be written as 
\begin{align}
\mathcal{Z}(o,s) \approx \mathcal{B}^{L}\mathcal{W}^{L}\mathcal{W}^{L-1}\ldots \mathcal{W}^1\mathcal{V}^0.
\label{eqn:butterfly_mat}
\end{align}
Here, the outer factor $\mathcal{V}^0=\mathrm{diag}(\mathcal{V}_{o,\nu_1},\ldots,\mathcal{V}_{o,\nu_{2^L}})$ consists interpolation matrices at the leafs of $\mathfrak{T}_s$, and the block-diagonal inner factors $\mathcal{W}^{l}, l=1,\ldots,L$ have blocks $\mathcal{W}_{\tau}$ for all nodes $\tau$ at level $l-1$ of $\mathfrak{T}_o$
\begin{align}
\mathcal{W}_\tau=
\begin{bmatrix}
\text{diag}(\mathcal{W}_{\tau_1,\nu_1}, \dots, \mathcal{W}_{\tau_1,\nu_{2^{L-l}}}) \\
\text{diag}(\mathcal{W}_{\tau_2,\nu_1}, \dots, \mathcal{W}_{\tau_2,\nu_{2^{L-l}}})
\end{bmatrix}
\end{align}
Here, $\nu_1, \nu_2, \ldots, \nu_{2^{L-l}}$ are the nodes at level $L-l$ of $\mathfrak{T}_s$ and $\tau_1$, $\tau_2$ are children of $\tau$. Lastly, the block-diagonal outer factor $\mathcal{B}^L=\mathrm{diag}(\mathcal{Z}(\tau_1,\bar{s}),\ldots,\mathcal{Z}(\tau_{2^L},\bar{s}))$ for all leafs $\tau_1, \tau_2, \ldots, \tau_{2^L}$ of $\mathfrak{T}_o$. Note that there exist several equivalent forms to \cref{eqn:butterfly_mat} \cite{liu2020butterfly,Yingzhou_2017_IBF,li2015butterfly,Pang2020IDBF}. Upon factorizing subblocks in $\mathcal{Z}$ in the form of \cref{eqn:butterfly_mat}, we obtain a compressed system matrix $\mathcal{Z}_{\chi_{con}}$. The storage and application costs of a matrix-vector product using $\mathcal{Z}_{\chi_{con}}$ scale as $O(n\log n)$. As a result, the overall HOD-BF construction requires $O(N\log^2 N)$ memory and CPU resources (\newrev{see Section III-C for more details}). The HOD-BF-compressed system matrix can then be used for fast matrix-vector multiplication during the iterative solution of the VIE.

\subsection{Inversion of HOD-BF-compressed system matrix \label{ssec:HODBFinversion}}
Once constructed, the inverse of the HOD-BF-compressed system matrix can be computed approximately and used as a preconditioner during the iterative solution of \cref{eq:4}. The inversion algorithm has been previously described in \cite{Liu_2017_HODBF,liu2020MF} and is briefly summarized as HODBF\_invert in \cref{alg:hodbf_invert}.

\newrev{For each non-leaf $\tau$, $\mathcal{D}_\tau=[\mathcal{D}_{\tau_1},\mathcal{B}_{\tau_1};\mathcal{B}_{\tau_2},\mathcal{D}_{\tau_2}]$. $\mathcal{D}_\tau^{-1}$ is computed by first recursively computing $\mathcal{D}_{\tau_1}^{-1}$ and $\mathcal{D}_{\tau_2}^{-1}$. Next, the butterfly block is updated as $\mathcal{B}_{\tau_i}\leftarrow \mathcal{D}^{-1}_{\tau_i}\mathcal{B}_{\tau_i}$ with both $\mathcal{D}^{-1}_{\tau_i}$ and $\mathcal{B}_{\tau_i}$ already compressed (see line 7). Finally the updated matrix $[\mathcal{I},\mathcal{B}_{\tau_1};\mathcal{B}_{\tau_2},\mathcal{I}]$ is inverted using the butterfly extension of the Sherman-Morrison-Woodbury formula \cite{Hager1989SMW}, named BF\_SMW in \cref{alg:bf_smw} (see line 8). }  

\newrev{The operations at line \ref{line:matvec1_hodbf} of \cref{alg:hodbf_invert} and lines \ref{line:matvec3_hodbf}, \ref{line:matvec4_hodbf} of \cref{alg:bf_smw}} require construction of a new butterfly representation out of combinations of existing butterfly representations. This type of operations can be accelerated by rapid multiplications of the operator with random vectors and reconstructing the butterfly from the results \newrev{with tolerance $\chi_{fact}$}. This algorithm is named BF\_random\_matvec and requires $O(n^{3/2}\log n)$ CPU operations \cite{liu2020butterfly}, and consequently the computational cost for $\mathcal{Z}^{-1}$ scales as $O(N^{3/2}\log N)$.

\begin{algorithm}
	\caption{${\rm HODBF\_invert}(\mathcal{A})$: Inversion of a HOD-BF matrix.}
	\label{alg:hodbf_invert}
	\hspace*{\algorithmicindent} \textbf{Input:} $\mathcal{A}$ in HOD-BF form \\
	\hspace*{\algorithmicindent} \textbf{Output:} $\mathcal{A}^{-1}$ in HOD-BF form 	
	\begin{algorithmic}[1]
		\State Let $\mathcal{D}_{\tau}=\mathcal{A}$ with $\tau$ denoting the root node.
		\If{$\mathcal{D}_{\tau}$ dense}
		\State Directly compute $\mathcal{D}_\tau^{-1}$
		\Else
		\State \newrev{$\mathcal{D}_\tau=[\mathcal{D}_{\tau_1},\mathcal{B}_{\tau_1};\mathcal{B}_{\tau_2},\mathcal{D}_{\tau_2}]$}
		\State $\mathcal{D}_{\tau_i}^{-1}\leftarrow$ HODBF\_invert$\left( \mathcal{D}_{\tau_i} \right)$, $i=1,2$
		\State $\mathcal{B}_{\tau_i}\leftarrow$ BF\_random\_matvec$\left( \mathcal{D}_{\tau_i}^{-1}\mathcal{B}_{\tau_i} \right)$, $i=1,2$\label{line:matvec1_hodbf}
		\State $\mathcal{D}_{\tau}^{-1}\leftarrow$ BF\_SMW$\left( \begin{bmatrix}
		\mathcal{I}& \mathcal{B}_{\tau_1}\label{line:smw_hodbf} \\
		\mathcal{B}_{\tau_2}&\mathcal{I} 
		\end{bmatrix} \right) $$\begin{bmatrix}\mathcal{D}_{\tau_1}^{-1}&\\&\mathcal{D}_{\tau_2}^{-1}\end{bmatrix}$
		\EndIf	
	\end{algorithmic}
	
\end{algorithm}

\begin{algorithm}
	\caption{${\rm BF\_SMW}(\mathcal{A})$: SMW inversion of a butterfly matrix.}
	\label{alg:bf_smw}
	\hspace*{\algorithmicindent} \textbf{Input:} $\mathcal{A}-\mathcal{I}$ is a butterfly of $L$ levels \\
	\hspace*{\algorithmicindent} \textbf{Output:} $\mathcal{A}^{-1}-\mathcal{I}$ is a butterfly of $L$ levels
	\begin{algorithmic}[1]
		\State Split $\mathcal{A}$ into four children butterflies of $L-2$ levels: $\mathcal{A}=[\mathcal{A}_{11},\mathcal{A}_{12};\mathcal{A}_{21},\mathcal{A}_{22}]$ using $\mathfrak{T}_H$\label{line:smw_split}
		\State $\mathcal{A}_{22}^{-1}\leftarrow$ BF\_SMW$\left( \mathcal{A}_{22} \right)$
		\State $\mathcal{A}_{11}\leftarrow$ BF\_random\_matvec$\left(\rev{ \mathcal{A}_{11}-\mathcal{A}_{12}\mathcal{A}_{22}^{-1}\mathcal{A}_{21}-\mathcal{I}}\right)$\label{line:matvec3_hodbf}
		\State $\mathcal{A}_{11}^{-1}\leftarrow$ BF\_SMW$\left(\mathcal{A}_{11}\right)$
		\State $\mathcal{A}^{-1}-\mathcal{I}\!\leftarrow$\! BF\_random\_matvec
		\Statex $ \left(\! \begin{bmatrix}
		\mathcal{I} & \\
		-\mathcal{A}_{22}^{-1}\mathcal{A}_{21} & \mathcal{I}
		\end{bmatrix}\!\!\!\begin{bmatrix}
		\mathcal{A}_{11}^{-1} & \\
		& \mathcal{A}_{22}^{-1}
		\end{bmatrix}\!\!\!\begin{bmatrix}
		\mathcal{I} &-\mathcal{A}_{12}\mathcal{A}_{22}^{-1} \\
		& \mathcal{I}
		\end{bmatrix}\!\!-\!\mathcal{I} \!\right)$\label{line:matvec4_hodbf}
	\end{algorithmic}
	
\end{algorithm}

\subsection{Rank scaling and complexity estimate\label{sec:CC}}

We first provide a simple analysis of the butterfly rank $r$ of the off-diagonal blocks $\mathcal{B}_\tau$ in HOD-BF, and then summarize the complexity of the proposed HOD-BF-accelerated solver with preconditioner.

The asymptotic scaling of the butterfly ranks has been studied for various schemes, including the strong-admissibility-based $\mathcal{H}$ matrix frameworks such as \cite{michielssen_multilevel_1996,Han_2017_butterflyLUPEC,Han_2018_butterflyLUDielectric}, $r=O(1)$ and weak-admissibility-based frameworks, e.g., HOD-BF. It was shown that $r=O(\log n)$ for 2-dimensional SIE \cite{Liu_2017_HODBF,Yang_2020_BFprecondition}, and $r=O(n^{1/4})$ for 3-dimensional SIE \cite{liu2020butterfly}. In what follows, we show that when applied to blocks in the 3-dimensional VIE system matrix, $r=O(n^{1/3})$. Moreover, despite of the non-constant butterfly ranks, the off-diagonal block can still be compressed with $O(n\log n)$ memory and CPU complexities using the butterfly algorithm.  

\begin{figure}[!htpb]
	\centering
	\includegraphics[width=0.95\columnwidth]{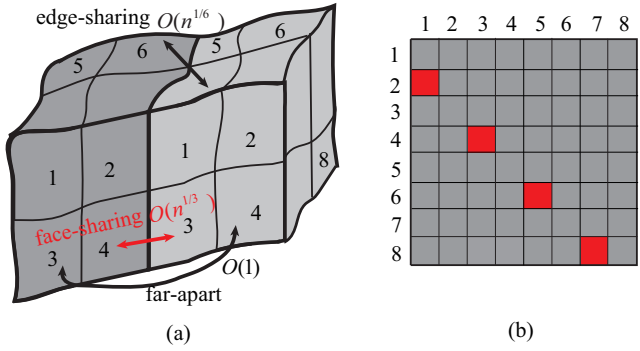}
	\caption[]{(a) Subscatterers (labeled from 1 to 8) at level $L/2$ of both $\mathfrak{T}_o$ and $\mathfrak{T}_s$ with $L=6$. (b) Subblocks at level $L/2$ of both $\mathfrak{T}_o$ and $\mathfrak{T}_s$. The subblocks with rank $O(n^{1/3})$ are highlighted in red. }
	\label{fig:rank}
\end{figure}

To begin with, \cref{fig:rank} shows an example of scatterer pair $(o,s)$ corresponding to one $m\times n$ block $\mathcal{Z}(o,s)$ with $L=6$, the subscatterers at level $L/2=3$ of $\mathfrak{T}_o$ and $\mathfrak{T}_s$ are shown in \cref{fig:rank}(a), and the subblocks corresponding to these subscatterer pairs are shown in \cref{fig:rank}(b). Apparently, there are $2^{L/2}=8$ subscatterers at level $L/2$ of both $\mathfrak{T}_o$ and $\mathfrak{T}_s$, and a total of $2^L=64$ subblocks $\mathcal{Z}(\tau,\nu)$ for nodes $\tau,\nu=1,\ldots,8$ at level $L/2$. Each subscatterer contains $O(n^{1/6})\times O(n^{1/6})\times O(n^{1/6})=O(n^{1/2})$ number of unknowns. Among the $2^L$ subblocks, the majority ($O(2^L)$ of them) representing interactions between well-separated subscatterer pairs have rank $r_{\tau,\nu}$ as $O(1)$ \cite{Bucci1987SpatialBandwidth,Engquist2018Greenfunction,Eric_1994_butterfly,Yang_2020_BFprecondition}. There are only $O(2^{L/3})\approx 4$ subblocks (highlighted in red) representing interactions between face-sharing subscatterer pairs. Note that here the big O notation is needed as in practise there can be pairs partially sharing one face (i.e., non-confirming partitioning in $\mathfrak{T}_o$ and $\mathfrak{T}_s$). The ranks of these subblocks are determined by the number of unknowns residing on the common face, leading to $r_{\tau,\nu}=O(n^{1/6})\times O(n^{1/6})=O(n^{1/3})$. Here, $O(n^{1/6})$ represents the number of unknowns in one of three dimensions of a subscatterer. In addition, there can be $O(2^{L/3})$ subblocks representing interactions between edge-sharing or vertex-sharing subscatterer pairs whose ranks scale as $O(n^{1/6})$ or $O(\log n)$. These edge-sharing or vertex-sharing pairs can be safely ignored compared to the face-sharing pairs. As a result, the transfer matrices $\mathcal{W}_{\tau,\nu}$ in $\mathcal{W}^{L/2}$ require only $O(2^L)\times (O(1))^2 + O(2^{L/3})\times (O(n^{1/3}))^2=O(n)$ storage units. Similarly, it can be shown that any factor $\mathcal{B}^L$, $\mathcal{W}^l$, and $\mathcal{V}^0$ requires $O(n)$ memory and hence the total memory for one butterfly is $O(n\log n)$. For the same reason, the butterfly construction of $\mathcal{Z}(o,s)$ requires $O(n\log n)$ CPU operations. \rev{BE CAREFUL: the entry evaluation is $O(r^2)$, the id is $O(r^3)$, but it seems entry evaluation has very large prefactors} Moreover, the randomized butterfly construction BF\_random\_matvec used in the inversion algorithm also requires $O(n^{3/2}\log n)$ CPU operations despite of the non-constant butterfly rank \cite{liu2020butterfly}. Following the analysis in \cite{Liu_2017_HODBF}, the computational costs of the proposed HOD-BF-accelerated solver and preconditioner for VIE scale as $O(N\log^2N)$ and $O(N^{3/2}\log N)$, respectively. 

\subsection{Parallelization\label{sec:parallel}}
\rev{
	The proposed distributed-memory parallelization scheme assigns processors according to the binary tree $\mathfrak{T}_H$. Let $\tau$ denote the root of $\mathfrak{T}_H$. Assuming a total of $P$ processors, the proposed scheme assigns all $P$ processors to $\tau$ and approximately $P/2$ processors to each of $\tau_1$ and $\tau_2$. Recall that $\mathcal{D}_\tau$ is partitioned into $[\mathcal{D}_{\tau_1},\mathcal{B}_{\tau_1};\mathcal{B}_{\tau_2},\mathcal{D}_{\tau_2}]$. The first $P/2$ processors handle $\mathcal{D}_{\tau_1},\mathcal{B}_{\tau_1}$, and the second $P/2$ processors handle $\mathcal{B}_{\tau_2},\mathcal{D}_{\tau_2}$. The processor subdivision continues until the subtree rooted at one node of $\mathfrak{T}_H$ is entirely handled by one processor. \newrev{Let $P_\tau$ denote the number of processors in charge of $\mathcal{D}_{\tau}$ and $\mathcal{B}_{\tau}$.} The butterfly represented off-diagonal blocks $\mathcal{B}_{\tau_1}, \mathcal{B}_{\tau_2}$ are distributed using the partitioning scheme detailed in \cite{liu2020butterfly}. \newrev{See \cref{fig:HOD-BF} for an example of the parallel layout with 8 processors.} According to this layout, Algorithms \ref{alg:hodbf_invert} and \ref{alg:bf_smw} are executed in parallel: 
	
	\newrev{
	In Algorithm 1, the update operation (line \ref{line:matvec1_hodbf}) are handled by $P_{\tau_i}$ processors. As $\mathcal{D}_{\tau_i}^{-1}$ and $\mathcal{B}_{\tau_i}$ are handled by the same set of processors, no data redistribution is needed. Since there is no data dependency between nodes at each level, Algorithm 1 is highly parallelizable. On the other hand, the BF\_SMW operation (line \ref{line:smw_hodbf}) is handled by $P_{\tau}$ processors. Due to the sequential nature of BF\_SMW, i.e., lines 2-5 of Algorithm 2 are executed sequentially, the parallelization scheme works as follows: assuming that BF\_SMW is executed on $P$ processors, a data redistribution is performed at the splitting step (line \ref{line:smw_split}) where each of $\mathcal{A}_{11}, \mathcal{A}_{12},\mathcal{A}_{21},\mathcal{A}_{22}$ is still distributed over $P$ processors following the layout of \cite{liu2020butterfly}. The splitting operation proceeds with $L/2$ levels until the children blocks become low-rank products. Note that during the splitting, the number of nonzero blocks in the children butterflies decreases and some of the $P$ processors become idle towards the leaf level of the splitting operation. To alleviate such load imbalance, each nonzero block in a butterfly can be distributed over multiple processors. Following the splitting step, no more redistribution of butterflies is needed at lines 2-5. In both Algorithm 1 and 2, the essential building block BF\_random\_matvec is parallelized following \cite{liu2020butterfly}. Note that for line 5 of Algorithm 2, the matrix-vector multiplication results are not conformal to the layout of $\mathcal{A}^{-1}$, hence a data redistribution of the vectors is also required. 
	
	Overall, this scheme allows good load balance and reduced collective communication during the overall HOD-BF construction, inversion and application. 
}
}

\section{Numerical Results}\label{sec:num_res}
This section presents several numerical examples that demonstrate the accuracy and efficiency of the proposed butterfly-accelerated VIE solver when applied to the EM analysis of heterogeneous scatterers.
Unless stated otherwise below, the scatterers are illuminated by an $\hat x$ polarized plane-wave travelling along $-\hat z$ direction.  All simulations are carried out on the Cori supercomputer of Haswell nodes; each node has two $16$-core Intel Xeon E5-2698v3 processors and $128$ GB of $2133$ MHz DDR4 memory. Unless otherwise stated, each simulation is carried out on 16 nodes.

\subsection{Validation}\label{sec:acc}
\subsubsection{Positive permittivity}
\begin{figure}[!ht]
\centering
\subfigure[]{\includegraphics[width=0.35\columnwidth]{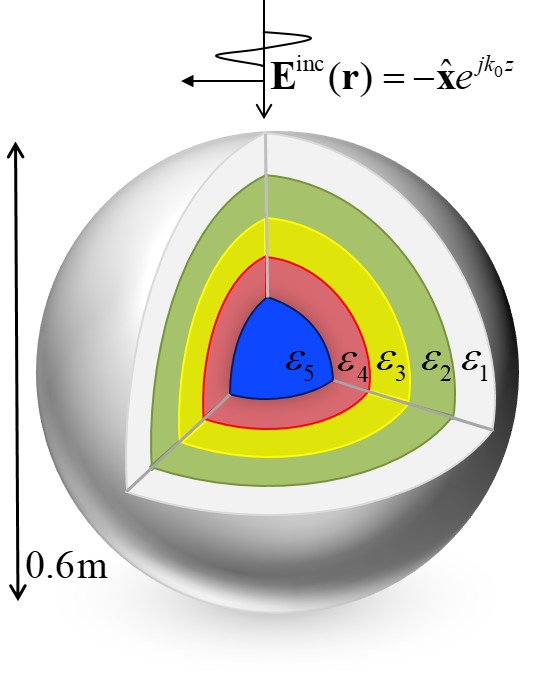}}
\subfigure[]{\includegraphics[width=0.58\columnwidth]{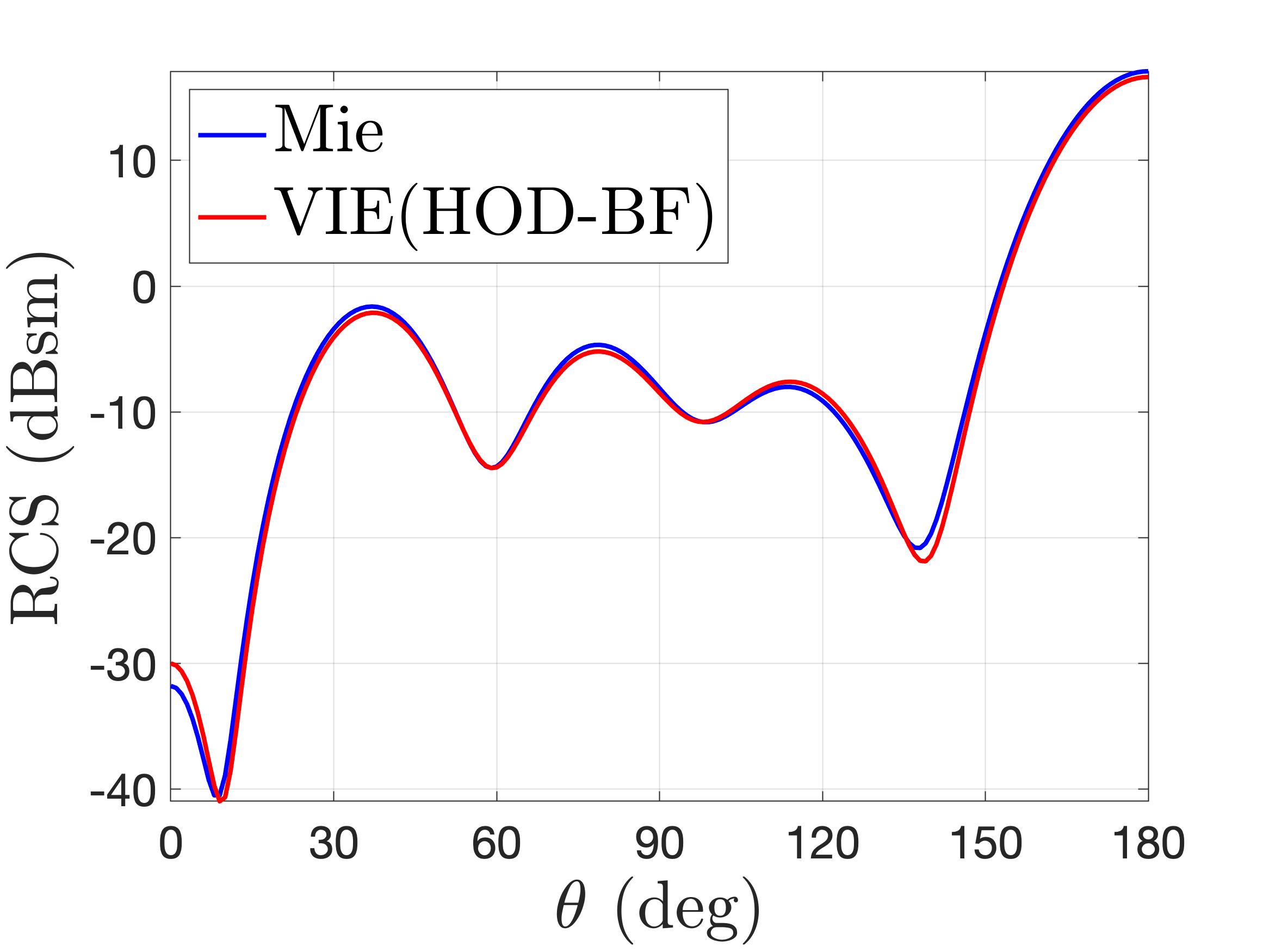}}
\caption[]{(a) The structure considered for the analysis: a five-layered sphere with $\varepsilon_1=2\varepsilon_0$, $\varepsilon_2=3\varepsilon_0$, $\varepsilon_3=4\varepsilon_0$, $\varepsilon_4=5\varepsilon_0$, and $\varepsilon_5=6\varepsilon_0$ . \rev{(b) The RCSs} computed by the proposed HOD-BF-accelerated VIE solver and Mie series solution at $f=750$ MHz.}
\label{fig:layered_sp_rcs}
\end{figure}

First, a five-layered sphere with layer thickness of $0.06$ m is analyzed. The permittivity of each layer, from outer to inner, is $2\varepsilon_0$, $3\varepsilon_0$, $4\varepsilon_0$, $5\varepsilon_0$ and $6\varepsilon_0$ [\cref{fig:layered_sp_rcs}(a)]. The sphere is discretised by a tetrahedral mesh with $N=407,842$, and the radar cross section (RCS) of the sphere is computed at $f=750$ MHz using the proposed VIE solver with compression tolerance $\chi_{con}=10^{-4}$, factorization tolerance $\chi_{fact}=10^{-1}$, and TFQMR stopping criteria $\chi_{sol}=10^{-3}$. The computed RCS is compared with the one obtained via (exact) Mie series solution [\cref{fig:layered_sp_rcs}(b)]. \rev{Let $x_i$ and $\hat{x}_i$, $i=1,\ldots,N_{rcs}$ denote the $N_{rcs}$ RCS samples computed by the proposed solver and Mie series solution, respectively. The relative root-mean-square error (RMSE) is defined as $\sqrt{\sum_i{(x_i-\hat{x}_i)^2}/N_{rcs}}/\max_i{\hat{x}_i}$, \newrev{where $x_i$ are in their original values rather than in dBsm}. The relative RMSE between two solutions is 0.024.} The good agreement between results shows that the proposed VIE solver with the chosen parameters achieve good accuracy. It should be noted that the proposed solver requires 2 and 78 TFQMR iterations with and without preconditioner, respectively. The details of the simulation performed by the proposed solver are provided in Table \ref{tab:tech_data}. 

\begin{figure}[!ht]
	\centering
	\subfigure[]{\includegraphics[width=0.35\columnwidth]{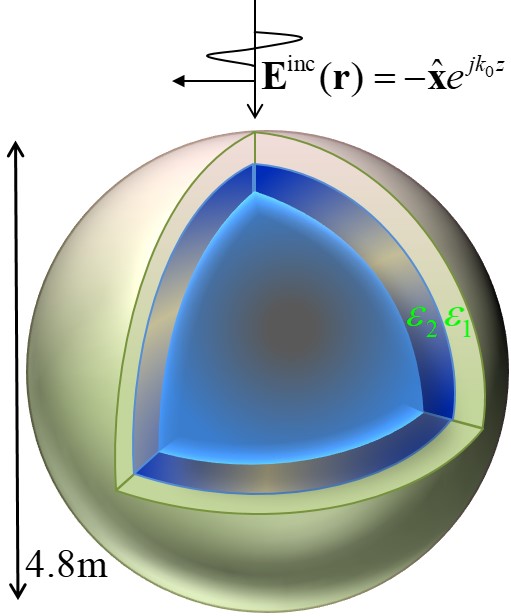}}
	\subfigure[]{\includegraphics[width=0.60\columnwidth]{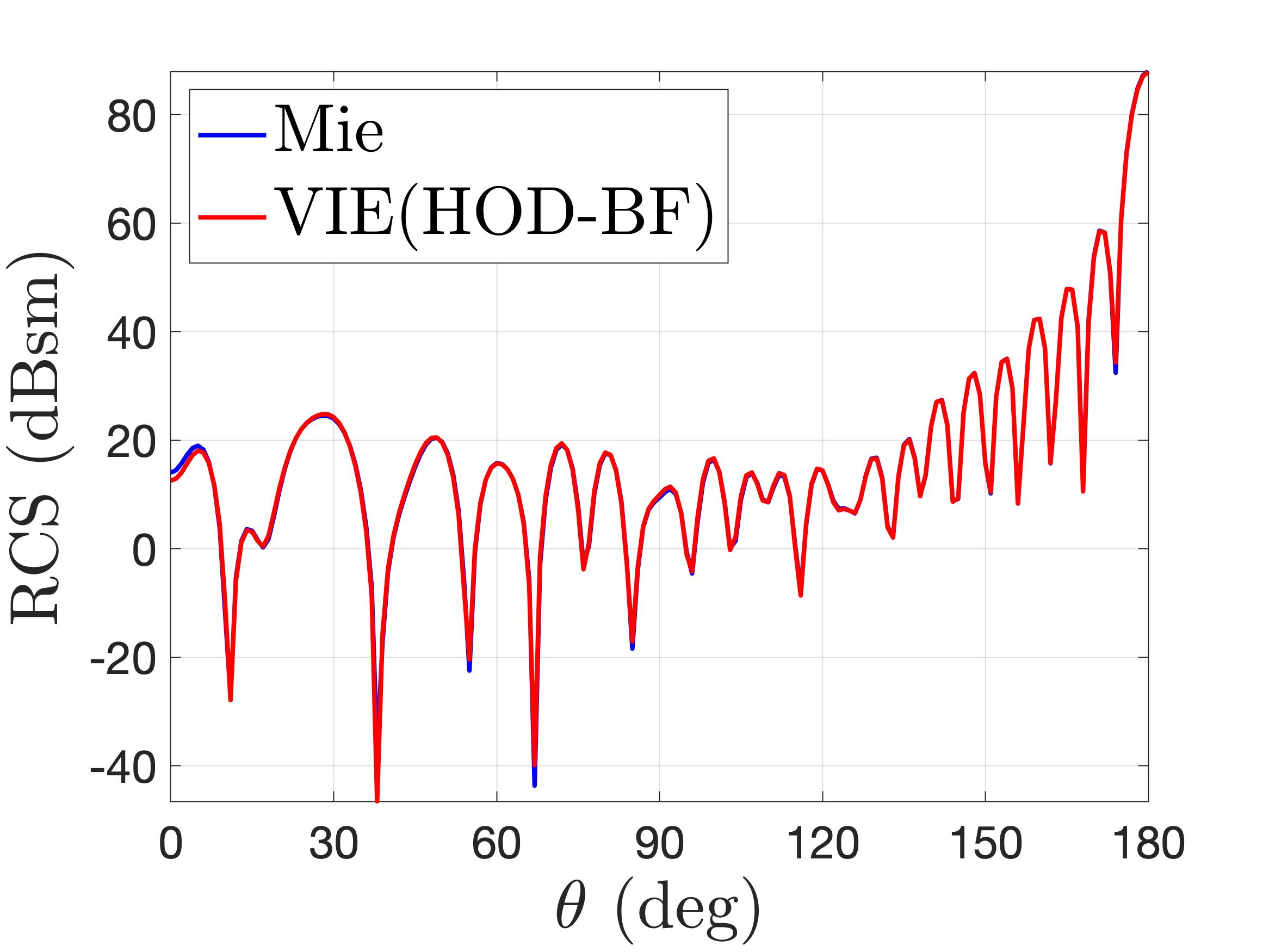}}
	\caption[]{(a)  The structure considered for the analysis: a two-layered shell with $\varepsilon_1=2\varepsilon_0$, and $\varepsilon_2=4\varepsilon_0$ . \rev{(b) The RCSs} computed with the proposed HOD-BF-accelerated VIE solver and Mie series solution at $f=600$ MHz.}
	\label{fig:2-layer_sp_shell_rcs}
\end{figure}

Next, a three-layered sphere with a hollow spherical core (i.e., two shells) is analyzed. The sphere has a radius of $2.4$ m and each shell has thickness of $0.03$ m. The dielectric permittivity of each shell, from inner to outer, is $4\varepsilon_0$ and $2\varepsilon_0$, respectively [\cref{fig:2-layer_sp_shell_rcs}]. The tetrahedral mesh of two shells consists of $N=5,530,950$ faces. The RCS is computed at $f=600$ MHz using the proposed VIE solver with  $\chi_{con}=10^{-4}$, $\chi_{fact}=10^{-1}$, and  $\chi_{sol}=10^{-3}$.  The computed RCS is compared with the one obtained via the Mie series solution [\cref{fig:2-layer_sp_shell_rcs}(b)]. \rev{The relative RMSE between two solutions is 0.0034.} Once again, the results show very good agreement. The proposed VIE solver with and without preconditioner requires 7 and 643 TFQMR iterations, respectively. More details on the simulation are provided in Table \ref{tab:tech_data}. 

	\begin{table}[htp!]
	\footnotesize
    \caption{Technical data for the numerical examples: layered sphere, shell, and sphere with negative permittivity.}
	\begin{center}
		\begin{tabular}{|c|c|c|c|}
			\hline			
			Example & Layered & Shell & Negative Sphere  \\
			\hline				
			Max. dimension & $1.5\lambda_0$ & $9.6\lambda_0$ & $1.6\lambda_0$ \\		
			\hline
			$N$ & $407,842$ & $5,530,950$ & $896,001$ \\			
			\hline						
			Max. rank & 1088 & 1331 & 1198  \\		
			\hline					
			Memory for $\mathcal{Z}$ & 57 GB & 324 GB & 96 GB  \\
			\hline							
			Memory for $\mathcal{Z}^{-1}$ & 28 GB & 321 GB & 306 GB  \\			
			\hline
			Construction time & 7.4 min & 1.5 h & 1.9 h   \\
			\hline								
			Inversion time & 26 min & 6.2 h & 6.5 h  \\
			\hline						
			Solve time w/ precond. & 2.3 s & {48.3} s  & 21.4 s  \\
			\hline						
			Iteration \# w/ precond. & 2 & 7 & 6  \\
			\hline						
			Solve time w/o precond. & {42.8 s} & {32} min & 46 min \\
			\hline						
			Iteration \# w/o precond. & {78} & 643  & 3447  \\						
			\hline
		\end{tabular}
	\end{center}
	\label{tab:tech_data}
\end{table}

\subsubsection{Negative permittivity}
\begin{figure}[!ht]
	\centering
\subfigure[]{\includegraphics[width=0.35\columnwidth]{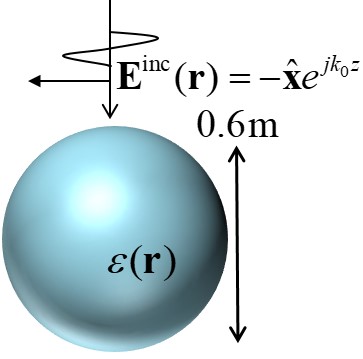}}
\subfigure[]{\includegraphics[width=0.60\columnwidth]{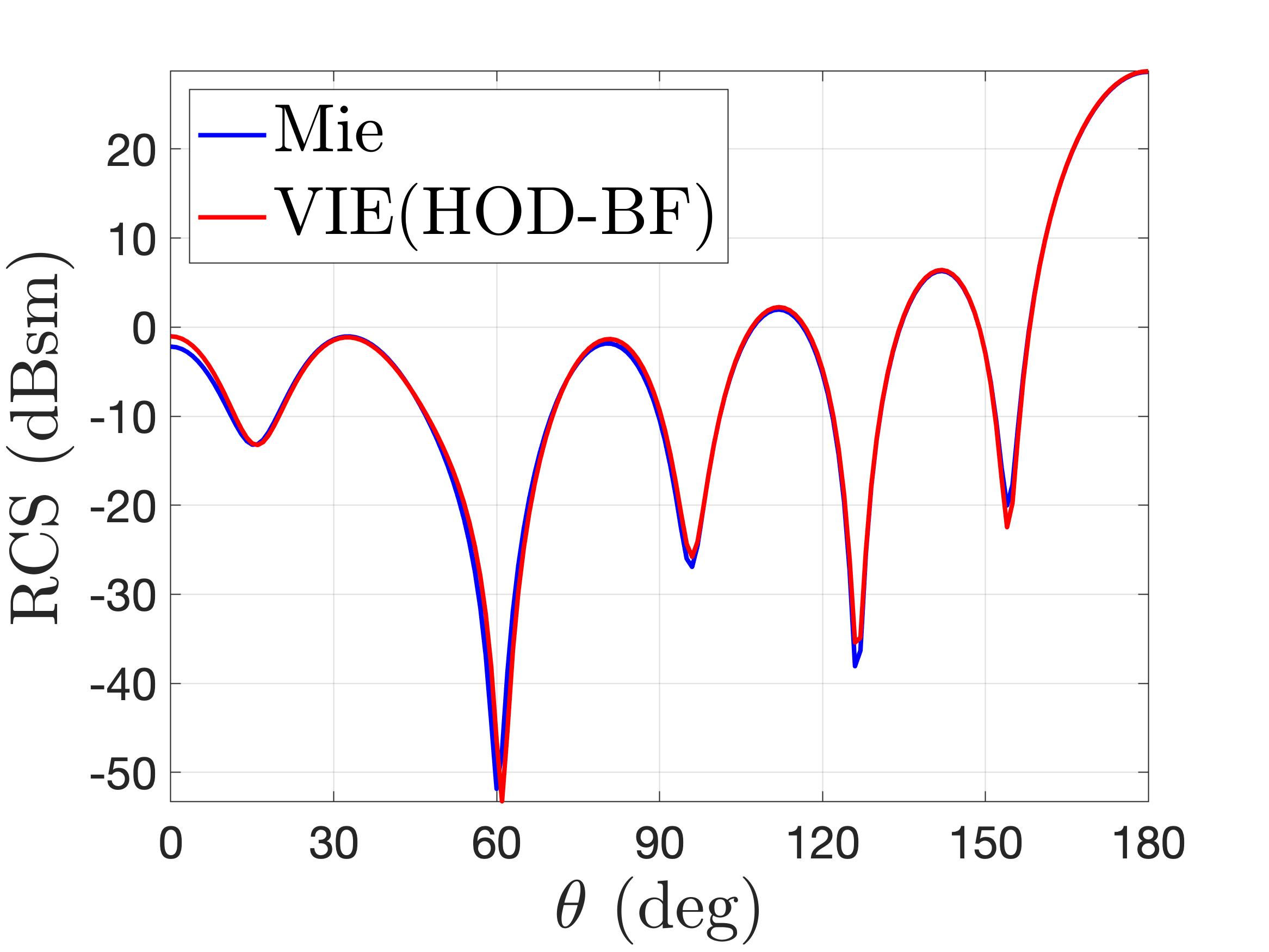}}
	\caption[]{(a) The structure considered for the analysis: a sphere with negative permittivity, $\varepsilon ({\mathbf{r}}) = \left( {-4 - 0.2j} \right){\varepsilon _0}$. \rev{(b) The RCSs} computed via the proposed HOD-BF-accelerated VIE solver and Mie series solution at  $f=800$ MHz.}
	\label{fig:neg_sp_rcs}
\end{figure}

\rev{ In general, the VIE formulation results in poorly conditioned systems when applied to the scattering problems involving negative permittivity. Iterative solution of such systems often do not converge with simple or no preconditioners. However, the use of HOD-BF with increased accuracy in factorization and solution provides stable and accurate results, even for the  ill-conditioned systems. To demonstrate the stability and accuracy of the proposed simulator for such systems,} a sphere of radius $0.3$ m with permittivity $\varepsilon ({\mathbf{r}}) = \left( { - 4 - 0.2j} \right){\varepsilon _0}$ is analyzed. 
The sphere is discretized by a tetrahedral mesh consisting of $N=896,001$ faces. The RCS is computed at $f=800$ MHz using the proposed VIE solver with  $\chi_{con}=10^{-4}$,  ${\chi_{fact}}=10^{-3}$, and $\chi_{sol}=10^{-5}$. The computed RCS is compared with the one obtained by the analytical Mie series solution [\cref{fig:neg_sp_rcs}(b)]; \rev{The relative RMSE between two solutions is 0.0027.} It is a well-known fact that the discretized VIE system becomes ill-conditioned for the EM analysis of the scatterers with negative permittivity and the iterative solution of the VIE system oftentimes cannot be obtained \cite{markkanen2016numerical,gomez2015volume,yucel2018internally}. Here the proposed VIE solver with preconditioner obtains the solution after 6 TFQMR iterations. \rev{Note that the choice of $\chi_{fact}=10^{-3}$ represents a trade-off between inversion time and solve time}. This clearly shows that the proposed solver is \rev{less sensitive} to the  ill-conditioning due to negative permittivity and it is robust for the board permittivity EM analysis. This example was run on 16 nodes and the details of the simulation are provided in Table \ref{tab:tech_data}.

\subsection{Computational Complexity}\label{sec:sca}
In this subsection, several examples were considered to verify the predicted CPU and memory scaling of the proposed HOD-BF-accelerated VIE solver. In all examples, the scatterers are discretized by the tetrahedal mesh with elements comparable to  ${\lambda _0}/10$. The performance of the proposed HOD-BF-accelerated VIE solver is compared with that of an $\mathcal{H}$ matrix-accelerated VIE solver with strong admissibility, and an \rev{HOD-LR-accelerated VIE solver}. \rev{For the $\mathcal{H}$ matrix-accelerated solver, all low-rank matrix blocks corresponding to $d_{a,b}>\alpha\max(r_a,r_b)$ are compressed, where $r_a$ and $r_b$ are the radii of the spheres enclosing subscatterers $a$ and $b$,  $d_{a,b}$ is the distance between the centers of subscatterers, and $\alpha=2$. We used a cobblestone distance sorting-based partitioning \cite{Shaeffer_2008} in the $\mathcal{H}$ matrix, HOD-LR, and HOD-BF accelerated solvers. For the low-rank compression algorithms, we used a blocked variant of adaptive cross approximation algorithm \cite{liu2019parallelACA} for better stability and parallel performance. The $\mathcal{H}$ matrix solver and HOD-LR solvers are available in the open-source software ButterflyPACK.} In the following tests,  ${\chi_{con}}=10^{-3}$,  ${\chi_{fact}}=10^{-1}$, and  $\chi_{sol}=10^{-3}$.

\subsubsection{Positive permittivity}

\begin{figure}[!ht]
\centering
\includegraphics[width=0.9\columnwidth]{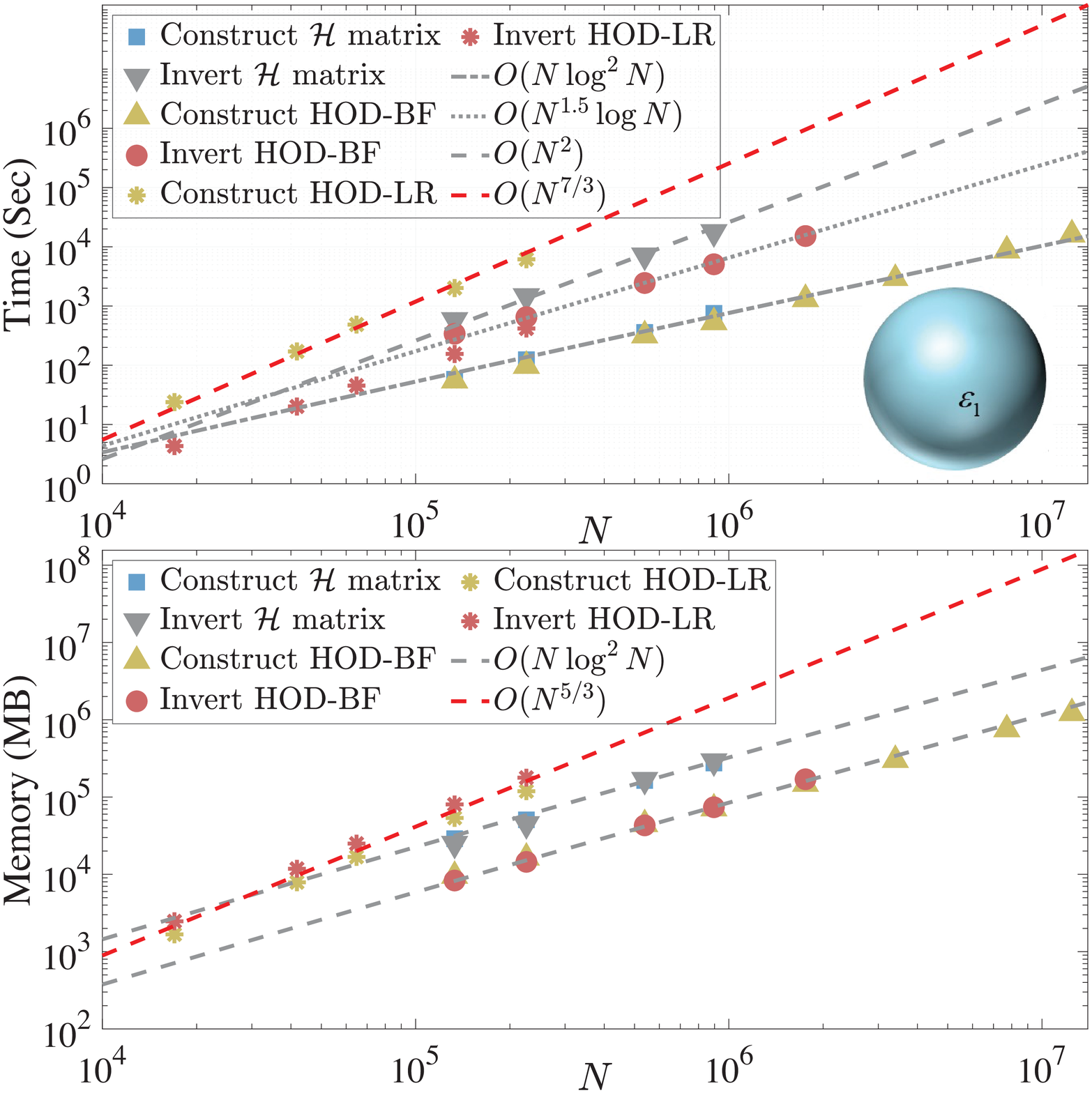} 
\caption[]{(a) \rev{The CPU and (b) memory scaling} of the HOD-BF-accelerated solver, $\mathcal{H}$ matrix-accelerated solver, \rev{and HOD-LR-accelerated solver} when applied to the EM analysis of a sphere with permittivity $\varepsilon ({\mathbf{r}}) = \left( {4 - 0.0001j} \right){\varepsilon _0}$.}
\label{fig:sp_scaling}
\end{figure}

\begin{figure}[!ht]
	\centering
	\includegraphics[width=0.9\columnwidth]{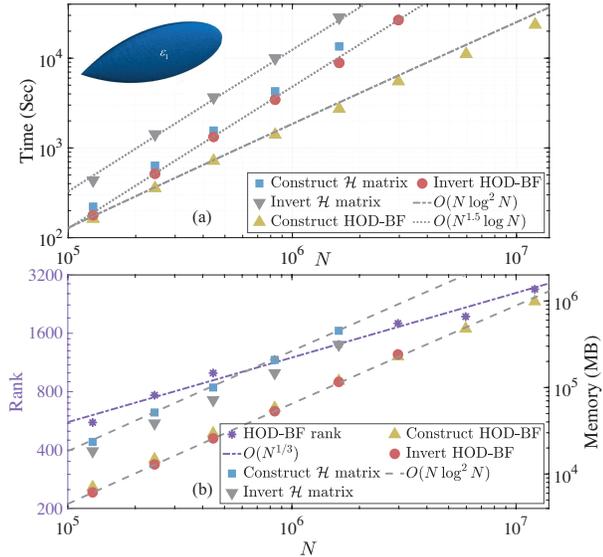} 
	\caption[]{(a) CPU and \rev{(b) memory scaling} of the HOD-BF-accelerated solver and  $\mathcal{H}$ matrix-accelerated solver when applied to the EM analysis of a NASA almond with permittivity $\varepsilon ({\mathbf{r}}) = \left( {4 - 0.0001j} \right){\varepsilon _0}$}
	\label{fig:nasa_scaling}
\end{figure}

\begin{figure}[!ht]
	\centering
	\includegraphics[width=0.9\columnwidth]{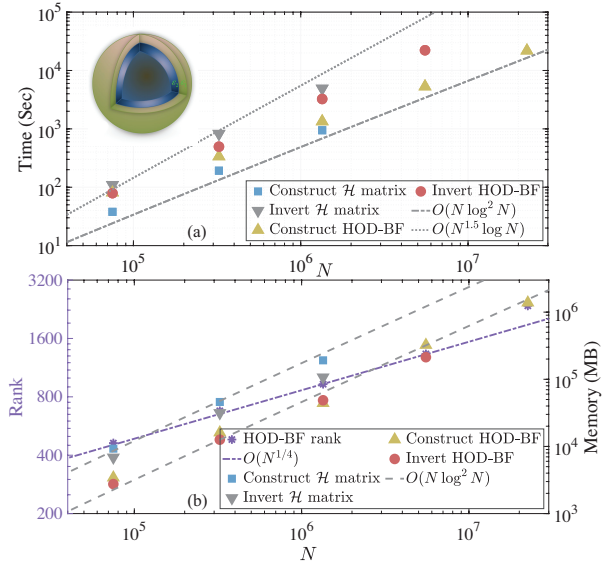} 
	\caption[]{(a) The CPU and \rev{(b) memory scaling} of the HOD-BF-accelerated solver and  $\mathcal{H}$ matrix-accelerated solver when applied to the EM analysis of a two-layered shell with inner and outer shell permittivities of  $4\varepsilon_0$ and $2\varepsilon_0$}
	\label{fig:shell_scaling}
\end{figure}

	\begin{table}[thp!]
	\footnotesize
    \caption{Time (in seconds) and iteration count for the solve phase of the proposed HOD-BF-accelerated VIE solver with (`w/') and without ('w/o') preconditioner. Note that each TFQMR iteration requires 6 MVPs.}
	\begin{center}
		\begin{tabular}{|c|c|c|c|c|}
			\hline			
			$N$ & $225,746$ & $538,782$ & $896,001$ & $1,757,024$\\
        	\hline						
			Solve time w/    & 1.3  & 3.3  & 6.77  & 17 \\
			\hline						
			Iteration \# w/  & 3 & 3 & 4 & 5 \\
			\hline
			Time per MVP w/  & 0.07  & 0.185  & 0.28   & 0.56 \\			
			\hline						
			Solve time w/o   &  48.9  &  310.1  & 509.2 & 1940 \\
			\hline						
			Iteration \# w/o   & 279 & 712  & 786  & 1535 \\	\hline
			Time per MVP w/o  & 0.03  & 0.07  & 0.11  & 0.21   \\	
			\hline
		\end{tabular}
	\end{center}
	
	\label{tab:tech_data_scaling_solve}
\end{table}

	\begin{table}[htp!]
	\footnotesize
    \caption{Technical data for the examples: sphere, almond, human head}
	\begin{center}
		\begin{tabular}{|c|c|c|c|}
			\hline			
			Example & Sphere & Almond & Head \\
			\hline				
			Max. dimension & $3.4\lambda_0$ & $10\lambda_0$ & $0.662\lambda_0$ \\		
			\hline
			$N$ & $1,757,024$ & $2,972,024$ & $1,759,208$\\			
			\hline						
			Max. rank & 1496 &  2006 &  1981   \\		
			\hline					
			Memory for $\mathcal{Z}$ & 148 GB & 228 GB & 314 GB  \\
			\hline							
			Memory for $\mathcal{Z}^{-1}$ & 170 GB &  241 GB &  353 GB  \\			
			\hline
			Construction time & 22 min &  1.5 h &  3 h  \\
			\hline								
			Inversion time & 4.2 h &  7.4 h &  7.3 h \\
			\hline						
			Solve time w/ precond. & 17 s & 24 s & 5.4 min \\
			\hline						
			Iteration \# w/ precond. & 5\label{8} & 4 & 90  \\
			\hline						
			Solve time w/o precond. & {32 min} &  {35 min} &  {1.8 h} \\
			\hline						
			Iteration \# w/o precond. & {1535} & {962} & {3698}   \\						
			\hline
		\end{tabular}
	\end{center}
	
	\label{tab:tech_data_scaling}
\end{table}

The computational complexity of the proposed scheme is studied as a function of electrical size by varying the frequency and fixing object size, or, by keeping the frequency fixed and varying the object size. In the first example, EM analysis of an homogeneous  sphere with radius  $0.3$ m  and permittivity $ \varepsilon ({\mathbf{r}}) = \left( {4 - 0.0001j}\right){\varepsilon _0}$ is considered for frequencies $700$ MHz, $850$ MHz, $1.2$ GHz, $1.3$ GHz, $1.6$ GHz, $2.0$ GHz, $2.5$ GHz, and $2.8$ GHz. The analysis at $2.8$ GHz requires a tetrahedral mesh with  ${N}=12,455,135$. Comparison of the complexities of the $\mathcal{H}$ matrix-accelerated solver and the proposed HOD-BF-accelerated solver  in  \cref{fig:sp_scaling} reveals that the CPU time for HOD-BF construction and inversion scale as $O(N\log^2 N)$ and $O(N^{1.5}\log N)$, respectively and the memory requirement scales as $O(N\log^2 N)$. These scaling results match with the predicted ones in \cref{sec:CC}. The maximum butterfly rank ranges from 745 to 2800, which scales as $O(N^{1/3})$. In contrast, the inversion and memory cost of $\mathcal{H}$ matrix-accelerated solver scale roughly $O(N^2)$ and $O(N\log^2 N)$, respectively. In addition, its construction cost scales with  $O(N^{1.5})$ , although the scaling plot was not provided in \cref{fig:sp_scaling} to avoid overcrowding.  In addition, the proposed HOD-BF-accelerated solver requires four times less memory compared to $\mathcal{H}$ matrix-accelerated solver. \rev{Moreover, the HOD-LR VIE solver scales as $O(N^{5/3})$ for memory and $O(N^{7/3})$ for CPU. Due to its prohibitively high CPU and memory requirement for 3D VIE, we do not consider HOD-LR solvers for the other experiments. }

The solution obtained via HOD-BF-accelerated solver without preconditioner, which only requires performing HOD-BF construction, will not converge to  $\chi_{sol}=10^{-3}$ within 3000 iterations when the electrical size of the object is large. For example, iterative solution via such solver requires 167 iterations for $N=13,3207$, 1535 iterations for $N=1,757,024$, and fails to converge for larger problems. In contrast, the HOD-BF-accelerated solver with preconditioner requires at most 5 iterations for the analysis at all frequencies. The data for the solve phase of HOD-BF-accelerated VIE solver with and without preconditioner is provided in \cref{tab:tech_data_scaling_solve}. As $N$ increases, the discretized VIE system becomes increasingly ill-conditioned, yet the approximate inverse preconditioner remains very effective and allows retaining the number TFQMR iterations a few. The time for each matrix-vector product (MVP) with the preconditioner is about twice of that without the preconditioner. The technical data for ${N}=1,757,024$ is provided in Table \ref{tab:tech_data_scaling}.

Next, a  NASA almond with size $0.25 \times 0.1 \times 0.04$ m and permittivity $ \varepsilon ({\mathbf{r}}) = \left( {4 - 0.0001j}\right){\varepsilon _0}$  is considered to study the computational complexity for frequencies $4.4$, $5.4$, $6.7$, $8$, $10$, $12$, $15$, and $20$ GHz. The analysis at $20$ GHz requires a tetrahedral mesh with ${N}=12,112,059$ . \cref{fig:nasa_scaling} shows that the CPU and memory requirements of the HOD-BF-accelerated solver for the construction and inversion scale as those predicted in \cref{sec:CC}. In contrast, the CPU scaling of the  $\mathcal{H}$ matrix-accelerated solver for the construction and inversion is at least $O(N^{1.5})$. Again, the iterative solution via HOD-BF-accelerated solver without preconditioner does not converge for the meshes with more than $N=2,972,024$ elements. However, the HOD-BF-accelerated solver with preconditioner requires maximum 4 TFQMR iterations for the analysis at all frequencies. More details on the analysis with the mesh comprising $N=2,972,024$ elements are provided in Table \ref{tab:tech_data_scaling}.

Furthermore, a three-layered sphere with a hollow spherical core (i.e., two shells) illuminated by a plane wave at $600$ MHz is considered for the complexity analysis. The permittivities of inner and outer shells with thickness of $0.03$ m are $4\varepsilon_0$ and $2\varepsilon_0$, respectively. The radius of the sphere is increased from $0.3$ m to $4.8$ m; the sphere with the largest radius requires a mesh with ${N}=22,491,763$. \cref{fig:shell_scaling} shows that the CPU and memory requirements of the HOD-BF-accelerated solver for the construction and inversion scale as those predicted in \cref{sec:CC}. \rev{Note that the memory requirement and construction time behave slightly better than the predicted $O(N\log^2N)$ scaling}. For this example, the scalings for the $\mathcal{H}$ matrix-accelerated solver are similar to those of the proposed solver, but with much larger leading constants. More details on the analysis with the mesh comprising $N=5,530,950$ elements are provided in  Table \ref{tab:tech_data}.

\subsubsection{Negative Permittivity}
In the next example, a sphere with permittivity $\varepsilon ({\mathbf{r}}) = \left( {-4 - 2j} \right){\varepsilon _0}$ is analyzed for frequencies $700$ MHz, $900$ MHz, $1.2$ GHz, $1.4$ GHz, and $1.7$ GHz. The analysis at $1.7$ GHz requires a mesh with  ${N}=1,757,024$ . 

\begin{figure}[!ht]
	\centering
	\includegraphics[width=0.95\columnwidth]{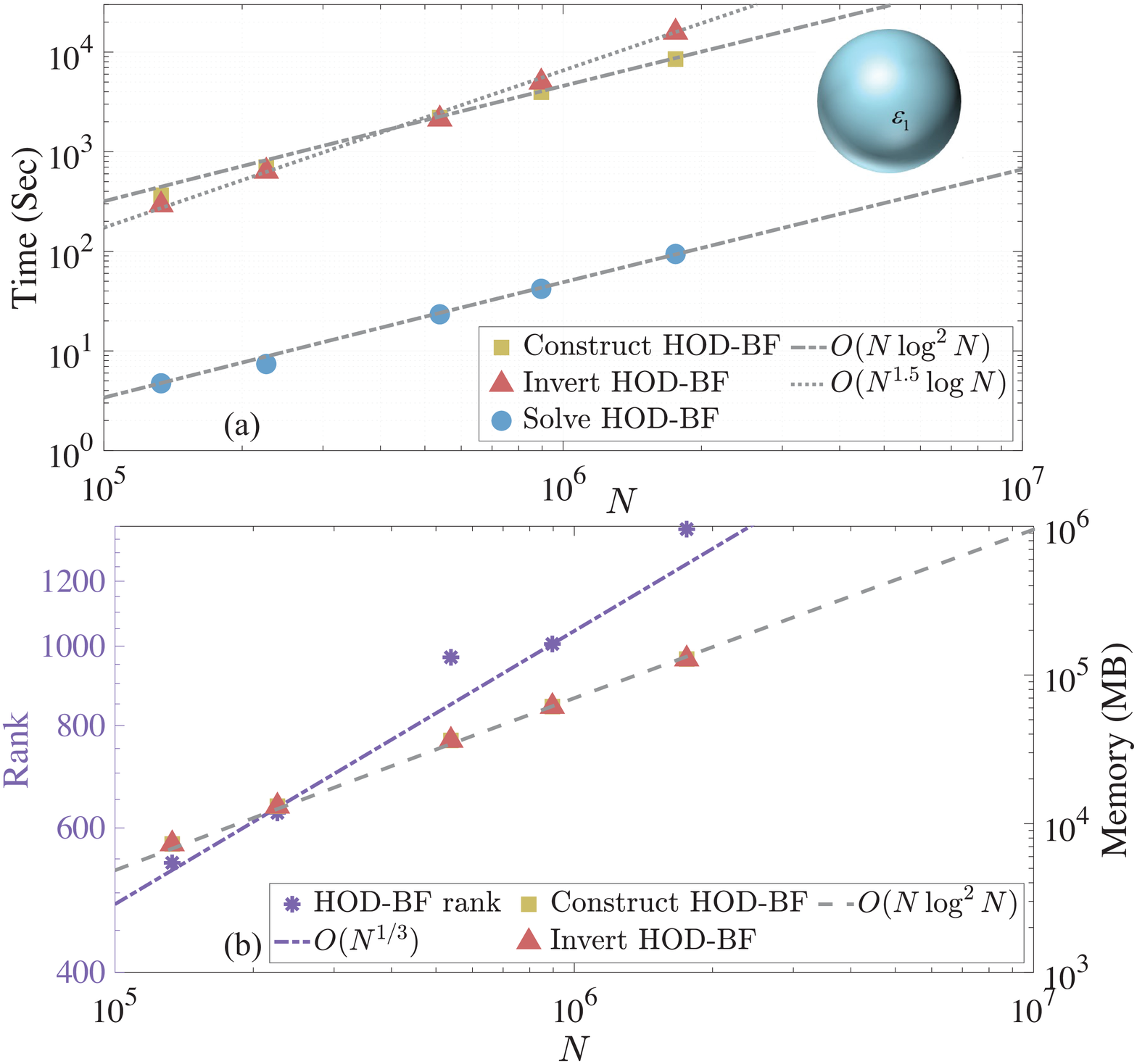} 
	\caption[]{(a) The CPU and \rev{(b) memory scaling} of the HOD-BF-accelerated solver and  $\mathcal{H}$ matrix-accelerated solver when applied to the EM analysis of  a sphere with permittivity $\varepsilon ({\mathbf{r}}) = \left( {-4 - 2j} \right){\varepsilon _0}$}
	\label{fig:neg_sp_scaling}
\end{figure}

\cref{fig:neg_sp_scaling} shows that the construction and inversion costs of the proposed HOD-BF-accelerated solver scale as $O(N\log^2 N)$ and $O(N^{1.5}\log N)$, respectively; its memory requirement scales as $O(N\log^2 N)$; the maximum butterfly rank scales as $O(N^{1/3})$. All the scaling results match the predicted ones in \cref{sec:CC}. Note that the solve time is also plotted (in blue) in \cref{fig:neg_sp_scaling}(a), which is negligible compared to the construction and inversion times, since the TFQMR iteration count is maximum 12 for the analyses at all frequencies. 

From the analysis in this subsection, we can conclude that the proposed HOD-BF-accelerated solver requires much less CPU time and memory compared to  $\mathcal{H}$ matrix-accelerated solver when the electrical size of the scatterer is large. Indeed, the complexities of $\mathcal{H}$ matrix-accelerated solver strongly depend on the geometry of the structure and can scale as bad as $O(N^2)$ . Furthermore, it is clear that the proposed simulator can provide efficient and accurate EM analysis even for the problems necessitating the solution of ill-conditioned VIE systems.

\subsection{Biomedical Application}
\begin{figure}[!ht]
	\centering
\subfigure[]{\includegraphics[width=0.49\columnwidth]{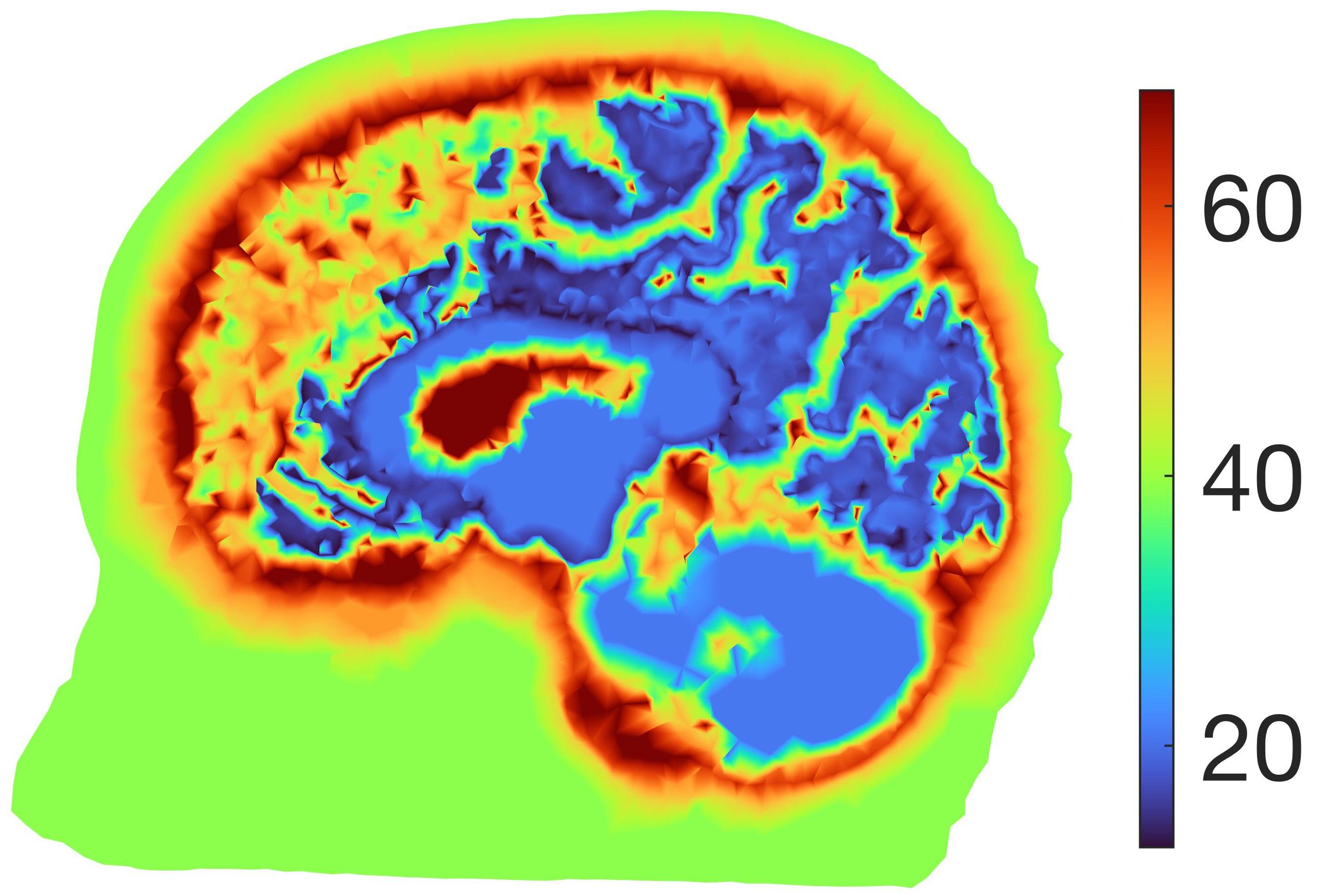}}
\subfigure[]{\includegraphics[width=0.49\columnwidth]{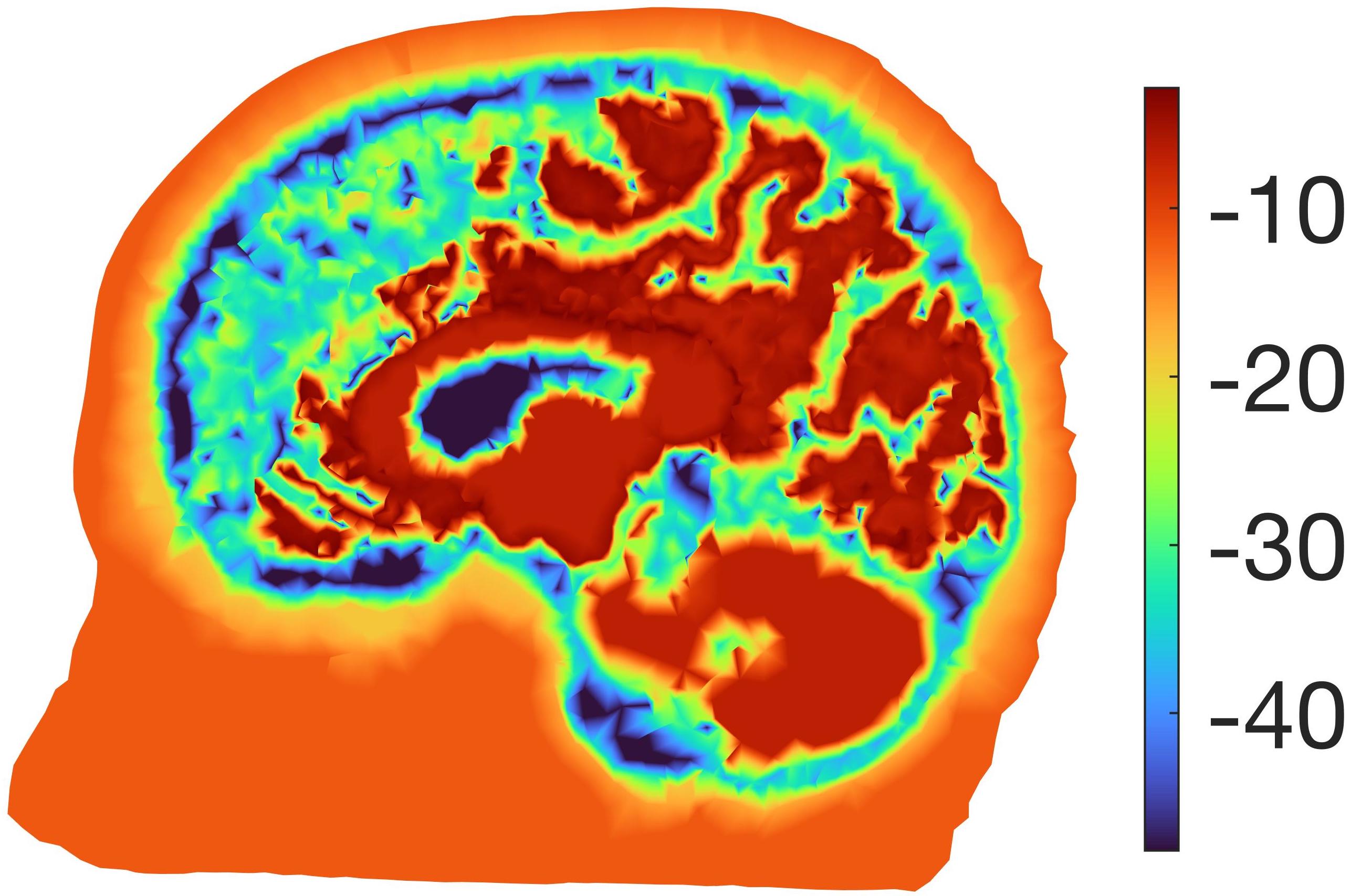}}
	\caption[]{The profile of the (a) real and (b) imaginary parts of the relative permittivities of tissues in the head model at $900$ MHz.}
	\label{fig:HeadExample}
\end{figure}

\begin{figure}[!ht]
	\centering
\subfigure[]{\includegraphics[width=0.49\columnwidth]{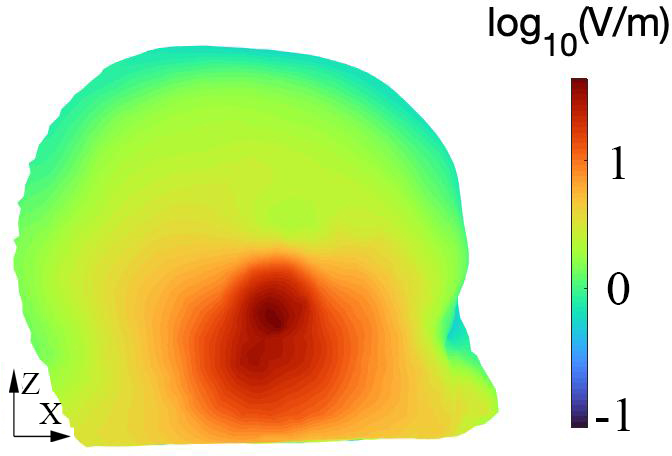}}
\subfigure[]{\includegraphics[width=0.46\columnwidth]{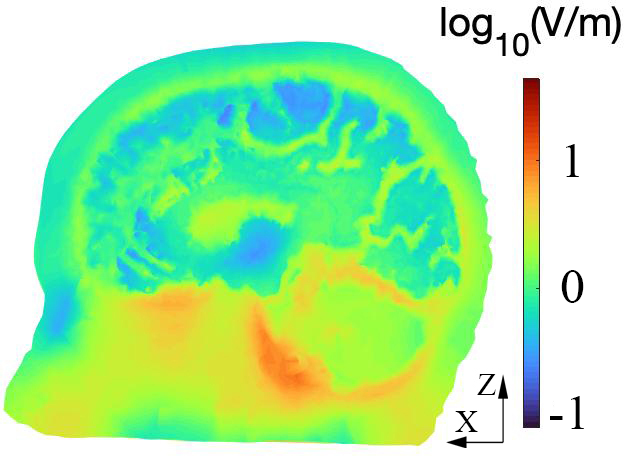}}
	\caption[]{\rev{(a) Right view and (b) left view of the magnitude of the electric field  (in log scale) inside the human head excited by a dipole with 1 \newrev{mA-m} strength at $900$ MHz computed by the proposed HOD-BF-accelerated solver.}}
\label{fig:FieldsInsideHead}
\end{figure}

The proposed HOD-BF-accelerated VIE solver is applied to the characterization of EM fields inside a human head induced by a cellphone operated at $900$ MHz. To this end, a head model is discretized using tetrahedrons with an average volume of 1 mm$^3$; the resulting mesh has  $N=1,759,208$ faces. The head model is excited by a \textit{z}-directed dipole positioned near to the ear. The relative permittivities of the tissues are plotted in Fig. \ref{fig:HeadExample}; those range from $12.50-j2.86$ to $68.60  -j48.134$.  The electric fields induced inside the human head is computed using the proposed HOD-BF-accelerated solver [Fig. \ref{fig:FieldsInsideHead}].

The HOD-BF-accelerated solver with $\chi_{con}=10^{-3}$, ${\chi_{fact}}=10^{-1}$ and $\chi_{sol}=10^{-5}$ requires 45 and 2568 TFQMR iterations with and without the preconditioner, respectively. This example was run on 32 nodes and the details of the simulation are provided in Table \ref{tab:tech_data_scaling}.

\section{Conclusion}\label{sec:con}
In this paper, a butterfly-accelerated VIE solver for the EM analysis of electrically-large scatterers comprising of broad permittivity values was proposed. The solver constructs the HOD-BF compressed blocks in the VIE system and accelerates the matrix-vector multiplications using these blocks. Furthermore, it leverages a preconditioner formed by the approximate inverse of the system matrix to ensure the fast convergence of the iterative solution. The complexity analysis and numerical experiments confirmed  $O(N\log^2N)$ CPU and memory scaling of the proposed solver for the construction of HOD-BF compressed blocks. Moreover, an $O(N^{1.5}\log N)$ CPU time scaling and $O(N\log^2N)$  memory scaling were validated for the inversion of the system matrix to form the preconditioner. The proposed HOD-BF-accelerated VIE solver was applied to EM analyses of various canonical and realistic scatterers. It was shown that the proposed solver can achieve significant CPU and memory savings compared to the low-rank-based $\mathcal{H}$ matrix-accelerated VIE solver for the high-frequency EM analysis.

\section*{Acknowledgment}\label{sec:ack}
Research reported in this publication was supported by the National Institute Of Mental Health of the National Institutes of Health under Award Number R00MH120046. The content is solely the responsibility of the authors and does not necessarily represent the official views of the National Institutes of Health. In addition, Y. Liu thanks the support from the U.S. Department of Energy, Office of Science, Office of Advanced Scientific Computing Research, Scientific Discovery through Advanced Computing (SciDAC) program through the
FASTMath Institute under Contract No. DE-AC02-05CH11231 at Lawrence Berkeley National Laboratory. 

%



\ifCLASSOPTIONcaptionsoff
  \newpage
\fi



\bibliographystyle{IEEEtran}
\bibliography{References.bib}
\end{document}